# **Instant Multiple Zeta Values at Non-Positive Integers and Bernoulli Functions** .

**Vivek V. Rane**

A-3/203 , Anand Nagar ,Dahisar ,

Mumbai-400 068

India

e-mail address : - v_v_rane@yahoo.co.in

**Abstract** : We instantly evaluate the multiple zeta function $\zeta_k(s_1, s_2, .........., s_k)$ at non-positive integers . We also deal with the Fourier theory of the products of Bernoulli polynomials . In particular, we obtain explicit generalized function B(s, $\alpha$ ) for the Bernoulli function . $B_m(\alpha)$ .

# Instant Multiple Zeta Values at Non-Positive Integers and Bernoulli Functions .

**Vivek V. Rane**

A-3/203 , Anand Nagar ,Dahisar ,

Mumbai-400 068

India

e-mail address : - v_v_rane@yahoo.co.in

For an integer $k \geq 1$ and for the complex variables $s_1, s_2, .........., s_k$ , let $\boxed{\zeta_k(s_1, s_2, .........., s_k)}$ be the multiple zeta function defined by $\boxed{\zeta_k(s_1, s_2, .........., s_k)}$

$= \sum_{n_1 \geq 1} n_1^{-s_1} \sum_{n_2 > n_1} n_2^{-s_2} ......... \sum_{n_k > n_{k-1}} n_k^{-s_k}$ for Re $s_1$, Re $s_2$,.........., Re $s_{k-1} \geq 1$ and Re $s_k > 1$ ;

and its analytic continuation . Then it is known that $\zeta_k(s_1, s_2, ........, s_k)$ is a meromorphic function of each of its variables . Our main object is to evaluate $\zeta_k(s_1, -m_2, -m_3, .........., -m_k)$, where $m_2, m_3, ..........m_k \geq 0$ are integers and $s_1$ is a complex number with Re $s_1 \leq 0$ ; or $s_1$ is a non-integral complex number with Re $s_1 > 0$ . Here we have defined $\zeta_k(s_1, -m_2, -m_3, ..........., -m_k)$

$= \lim_{s_1' \to s_1} \lim_{s_2' \to -m_2} .......... \lim_{s_k \to -m_k} \zeta_k(s_1', s_2', ..........s_k')$ , where $s_1' s_2', ..........s_k'$ are complex variables .

For a fixed complex number $\alpha \neq 0, -1, -2, -3, .......$ and for a complex variable $s$ , we define Hurwitz zeta function $\zeta(s, \alpha)$ by $\zeta(s, \alpha) = \sum_{n \geq 0} (n + \alpha)^{-s}$ for Re $s > 1$ ; and its analytic continuation . We write $\zeta(s,1) = \zeta(s)$ , the Riemann zeta function . Then it is known that $\zeta(s,\alpha)$ is a meromorphic function of $s$ with a simple pole at $s=1$ with residue 1 . From author [3] , [4] is clear that $\zeta(s,\alpha)-\zeta(s)$ is an analytic function of $\alpha$ except possibly for $\alpha \neq 0,-1,-2,..........$ More generally from author [3] , $\zeta^{(r)}(s,\alpha)-\zeta^{(r)}(s)$ is an analytic function of $\alpha$ for a fixed $s$ , if $\alpha \neq 0,-1,-2,..........$, where $\zeta^{(r)}(s,\alpha)=\frac{\partial^r}{\partial s^r}\zeta(s,\alpha)$. However $\zeta(-m,\alpha)$ for integral $m \geq 0$, is an entire function of $\alpha$ . We have the following three facts from author [3] , concerning $\zeta(s,\alpha)$, which we state as Lemmas .

**Lemma 1** : We have for any integer $m \geq 0$ and for $0 \leq \alpha \leq 1$ ,

$$\zeta(-m,\alpha)=m!\sum_{|n|\geq 1}\frac{e^{2\pi i n\alpha}}{(2\pi i n)^{m+1}} \quad .$$

**Note** : this is a consequence of Lemma 4 of author [3] .

**Lemma 2** : We have for any integer $m \geq 0$, $\zeta(-m,\alpha)-\zeta(-m,\alpha+1)=\alpha^m$ .

**Note** : From author [3] , this follows from $\zeta(s,\alpha)-\zeta(s,\alpha+1)=\alpha^{-s}$ for $\alpha \neq 0,-1,-2.-3,..........$ and for any complex $s$ .

**Lemma 3** : We have for any integer $m \geq 0$ and for any complex number $\alpha$ ,

$$\zeta(-m,\alpha)=\sum_{k=0}^{m}\binom{m}{k}\zeta(-k)\alpha^{m-k}+\alpha^m-\tfrac{\alpha^{m+1}}{m+1} \text{ , where } \binom{0}{0}=1 \text{ .}$$

**Note** : This is Theorem 1 of author [3] .

Next , the Bernoulli polynomial of degree $n \geq 0$ in variable $\alpha$ is defined by the expression $\frac{ze^{\alpha z}}{e^z-1}=\sum_{n\geq 0}\frac{B_n(\alpha)}{n!}z^n$ . for $|z|<2\pi$. Note that $B_0(\alpha)=1$. The Bernoulli numbers $B_n$ are defined by $B_n = B_n(0)$ . We shall write $\overline{B}_n = B_n(1)$. Then it is clear from Lemma 2* below that $\overline{B}_n = B_n$ for $n \neq 1$ and $\overline{B}_1 = B_1 + 1$ . We shall see that for any integer $m \geq 0$, $\zeta(-m, \alpha+1) = -\frac{1}{m+1}\sum_{k=0}^{m+1}\binom{m+1}{k}\overline{B}_k \alpha^{m+1-k}$ . This fact will be used in the proof of our Theorem .

We have following three facts concerning $B_m(\alpha)$, which we state as Lemmas .

**Lemma 1*** : We have for any integer $m \geq 0$ and for $0 \leq \alpha \leq 1$ ,

$$B_m(\alpha) = -(m+1)!\sum_{|n|\geq 1}\frac{e^{2\pi i n\alpha}}{(2\pi i n)^{m+1}} .$$

**Lemma 2*** : We have for any integer $m \geq 0$ , $B_m(\alpha+1) - B_m(\alpha) = m\alpha^{m-1}$ .

**Lemma 3*** : We have for any integer $m \geq 0$, $B_m(\alpha) = \sum_{k=0}^{m}\binom{m}{k}B_k \alpha^{m-k}$ .

From Lemma 1 and Lemma 1* , we have $\zeta(-m, \alpha) = -\frac{B_{m+1}(\alpha)}{m+1}$ for any integer $m \geq 0$ so that $\zeta(-m) = \zeta(-m, 1) = -\frac{\overline{B}_{m+1}}{m+1}$ . Hence , Lemma 2 and Lemma 2* are equivalent . In author [3] , we have stated that Lemma 3 and Lemma 3* are equivalent . This will be proved as Proposition 1 in what follows . In what follows , $\Gamma(s)$ will stand for gamma function .

For a complex $\alpha \neq 0, -1, -2, \ldots\ldots\ldots$, consider the function $B(s, \alpha) = s\zeta(s+1, \alpha)$. Then as a function of the complex variable $s$ , $B(s, \alpha)$ is an entire function of $s$ with $B(0, \alpha) = \lim_{s\to 0} s\zeta(s+1, \alpha) = 1$ . The fact that this limit is equal to 1 , can be proved by taking $0 < \alpha \leq 1$ to begin with ; and then extending the result for any complex $\alpha$ . Thus

$B(0,\alpha)$ is the constant polynomial 1 . If $m \ge 0$ is an integer , then

$B(-m,\alpha) = -m\zeta(-(m-1),\alpha) = B_m(\alpha)$, which is a consequence of Lemma 1 and Lemma 1* . Note that this is valid even for $m = 0$. Thus $B(s,\alpha)$ is the natural extension of the Bernoulli polynomial $B_m(\alpha) = B(-m,\alpha)$ , where $m \ge 0$ is an integer .

Akiyama and Tanigawa [1] have evaluated $\zeta_k(0,0,.........,0,-m)$ and $\zeta_k(-m,0,0,.........,0)$, where $m \ge 0$ is an integer . In this context , we refer the reader to Akiyama , Egami , Tanigawa [2] also .

Next , we state our Theorem . Though our result is true for any integer $k \ge 1$ , we shall state and prove our Theorem in the case of $k = 3$ for the sake of simplicity and convenience .

**Theorem** : Let $s_1$ be a complex number such that Re $s_1 \le 0$ ; or if Re $s_1 > 0$ , let $s_1$ non-integral . Let $m_2, m_3 \ge 0$ be integers . For any integer $k \ge 0$, let $\overline{B}_k = B_k(1)$ , where $B_k(\alpha)$ is the Bernoulli polynomial of degree k .Then

$$\zeta_3(s_1,-m_2,-m_3) = \sum_{k_3=0}^{m_3+1} \binom{m_3+1}{k_3} \frac{\overline{B}_{k_3}}{m_3+1} \cdot \sum_{k_2=0}^{m_2+m_3+2-k_3} \binom{m_2+m_3+2-k_3}{k_2} \frac{\overline{B}_{k_2}}{m_2+m_3+2-k_3} \cdot$$
$$\cdot \zeta(s_1-(m_2+m_3+2-k_2-k_3))$$

**Proposition 1** : For any integer $m \ge 0$ and for any complex number $\alpha$ , the following two statements are equivalent .

I) $\zeta(-m,\alpha) = \sum_{k=0}^{m} \binom{m}{k} \zeta(-k)\alpha^{m-k} + \alpha^m - \frac{\alpha^{m+1}}{m+1}$

II) $B_{m'}(\alpha)=\sum_{k=0}^{m'}\binom{m'}{k}B_k\alpha^{m'-k}$ , where $m'=m+1$ .

Our proof of Theorem will depend upon the fact that

$$\zeta(-m,\alpha+1)=-\tfrac{B_{m+1}(\alpha+1)}{m+1}$$

$$=-\tfrac{1}{m+1}\left\{B_{m+1}(\alpha)+(m+1)\alpha^m\right\}=-\tfrac{1}{m+1}\left\{\sum_{k=0}^{m+1}\binom{m+1}{k}B_k\alpha^{m+1-k}+(m+1)\alpha^m\right\}$$

$=-\tfrac{1}{m+1}\sum_{k=0}^{m+1}\binom{m+1}{k}\overline{B}_k\alpha^{m+1-k}$ , in view of the fact that $B_k=\overline{B}_k$ for $k\neq 1$ and $\overline{B}_1=1+B_1$ .

As $B_{m+1}(\alpha)=\sum_{k=0}^{m+1}\binom{m+1}{k}B_k\alpha^{m+1-k}$ , we shall write $\overline{B}_{m+1}(\alpha)$ for $\sum_{k=0}^{m+1}\binom{m+1}{k}\overline{B}_k\alpha^{m+1-k}$

Thus we have $\zeta(-m,\alpha+1)=-\frac{\overline{B}_{m+1}(\alpha)}{m+1}$ for any integer $m\geq 0$ .

Next , we shall discuss the Fourier theory of $\zeta(s,\alpha)$ for Re $s<1$ , as a function of $\alpha$ on the unit interval [0,1] . This has been already discussed in Lemma 4 of author [3] and in author [4] . As has been noted there ,

$\zeta(s,\alpha)\approx\Gamma(1-s)\sum_{|n|\geq 1}e^{2\pi in\alpha}(2\pi in)^{s-1}$ is the Fourier series of $\zeta(s,\alpha)$ as a function of $\alpha$ on the unit interval [0,1] , when Re $s<1$ . Here $\Gamma$ denotes gamma function . This is actually

Hurwitz's formula for $\zeta(s,\alpha)$ written in a suitable form . As the right hand side series in $\alpha$ converges boundedly to $\zeta(s,\alpha)$ on the interval [0,1] for Re $s<1$ , the right hand series is actually the Fourier series of $\zeta(s,\alpha)$ for Re $s<1$ . Using Parseval's Theorem for the product of two functions , we have , if Re $s_1$, Re $s_2<1$ , then

$$\int_0^1 \zeta(s_1,\alpha)\zeta(s_2,\alpha)d\alpha = (2\pi)^{s_1+s_2-2}\Gamma(1-s_1)\Gamma(1-s_2)\sum_{|n_1|,}\sum_{|n_2|\geq 1}\cdot\int_0^1 e^{2\pi i(n_1+n_2)\alpha}(in_1)^{s_1-1}(in_2)^{s_2-1}\,\mathrm{d}\alpha$$

$$= \Gamma(1-s_1)\Gamma(1-s_2)\sum_{\substack{|n_1|,\\ n_1=-n_2}}\sum_{|n_2|\geq 1}(2\pi i n_1)^{s_1-1}(2\pi i n_2)^{s_2-1} = 2(2\pi)^{s_1+s_2-2}\cos\tfrac{\pi}{2}(s_1-s_2)\Gamma(1-s_1)\Gamma(1-s_2)\zeta(2-(s_1+s_2))$$

if Re $(s_1+s_2)<1$ .

From Lemma 1* , we have $B_{m+1}(\alpha) = -(m+1)!\sum_{|n|\geq 1}\frac{e^{2\pi i n\alpha}}{(2\pi i n)^{m+1}}$ .

More generally , if $m_1, m_2, .....m_r \geq 0$ are integers , then the product

$B_{m+1}(\alpha)\ B_{m+2}(\alpha).......B_{m_r+1}(\alpha)$ has the Fourier series on the interval [0,1] namely ,

$$(-1)^r\prod_{i=1}^{r}(m_i+1)\,!\cdot\sum_{|n_1|\geq 1}\frac{e^{2\pi i n_1\alpha}}{(2\pi i n_1)^{m_1+1}}\sum_{|n_2|\geq 1}\frac{e^{2\pi i n_2\alpha}}{(2\pi i n_2)^{m_2+1}}..........\sum_{|n_r|\geq 1}\frac{e^{2\pi i n_r\alpha}}{(2\pi i n_r)^{m_r+1}}$$

$$=(-1)^r\prod_{i=1}^{r}(m_i+1)\,!\sum_{|n_1|,}\sum_{|n_2|}.........\sum_{|n_r|\geq 1}\left(\prod_{i=1}^{r}(2\pi i n_i)^{-(m_i+1)}\right)e^{2\pi i(\sum_i n_i)\alpha} = \sum_{|N|\geq 0}a_N e^{2\pi i N\alpha},\ \text{say}.$$

Thus , $\prod_{i=1}^{r}B_{m_2+1}(\alpha)$ has $\sum_{|N|\geq 0}a_N e^{2\pi i N\alpha}$ as the Fourier series on the interval [0,1] .

Next , $\int_0^1\prod_{i=1}^{r}B_{m_2+1}(\alpha)d\alpha = (-1)^r\prod_{i=1}^{r}(m_i+1)!\sum_{|n_1|,}\sum_{\substack{|n_2|,\\ \sum_i n_i=0}}......\sum_{|n_r|\geq 1}(\prod_{i=1}^{r}(2\pi i n_i)^{-(m_i+1)})$

On noting $B_{m_i+1}(\alpha) = \sum_{k_i=0}^{m_i+1}\binom{m_i+1}{k_i}B_{k_i}\alpha^{m_i+1-k_i}$ , and on integrating the polynomial expression

for $\prod_{i=1}^{r}B_{m_i+1}(\alpha)$ with respect to $\alpha$ from 0 to 1, we have proved the following

**<u>Proposition</u> 2** : We have $\sum_{k_1=0}^{m_1+1}\cdot\sum_{k_2=0}^{m_2+1}.........\sum_{k_r=0}^{m_r+1}\frac{\prod_{i-1}^{r}B_{k_i}\binom{m_i+1}{k_i}}{(1+\sum_i\ (m_i+1-k_i))}$

$$=(-1)^r \prod_{i=1}^{r} (m_i+1)! \sum_{\substack{|n_1|\ge 1 \\ \sum_{i=1}^{r} n_i = 0}} (2\pi i n_1)^{-(m_1+1)} \sum_{|n_2|\ge 1} (2\pi i n_2)^{-(m_2+1)} \ldots\ldots\ldots \sum_{|n_r|\ge 1} (2\pi i n_r)^{-(m_r+1)}$$

Next , we prove Proposition 1 .

**<u>Proof</u> <u>of</u> <u>Proposition</u> 1 :**

We shall show that statement I) implies statement II) .

$\zeta(-m,\alpha) = -\frac{B_{m+1}(\alpha)}{m+1}$ for $0<\alpha\le 1$ and for $m \ge 0$ . In particular , this means

$\zeta(-k) = -\frac{\overline{B}_{k+1}}{k+1}$ for $k \ge 0$ , where $\overline{B}_k = B_k(1)$ . Next , we have for $m \ge 0$ ,

$$\zeta(-m,\alpha) = \sum_{k'=0}^{m} \binom{m}{k'} \zeta(-k')\alpha^{m-k'} + \alpha^m - \frac{\alpha^{m+1}}{m+1} = -\frac{1}{m+1}\left\{ \alpha^{m+1} - (m+1)\alpha^m - (m+1)\sum_{k'=0}^{m}\binom{m}{k'}\zeta(-k')\alpha^{m-k'} \right\}$$

$$= -\frac{1}{m+1}\left\{ \alpha^{m+1} - (m+1)\alpha^m + (m+1)\sum_{k'=0}^{m}\binom{m}{k'}\frac{\overline{B}_{k'+1}}{k'+1}\alpha^{m-k'} \right\}$$

$$= -\frac{1}{m+1}\left\{ \alpha^{m+1} - (m+1)\alpha^m + \sum_{k'=0}^{m}\binom{m+1}{k'+1}\overline{B}_{k'+1}\alpha^{(m+1)-(k'+1)} \right\} .$$

Writing $k'+1 = k$ ,we have

$$\zeta(-m,\alpha) = -\frac{1}{m+1}\left\{ \alpha^{m+1} - (m+1)\alpha^m + \sum_{k=1}^{m+1}\binom{m+1}{k}\overline{B}_k\alpha^{m+1-k} \right\}$$

$$= -\frac{1}{m+1}\left\{ \sum_{k=0}^{m+1}\binom{m+1}{k}\overline{B}_k\ \alpha^{m+1-k} - (m+1)\alpha^m \right\} .$$

From Lemma 2* , we have $B_n(1) = B_n(0)$ for $n \ne 1$ and $B_1(1) = B_1(0)+1$ .

Thus we have $\zeta(-m,\alpha) = -\frac{1}{m+1}\left\{ \sum_{\substack{k=0 \\ k\ne 1}}^{m+1}\binom{m+1}{k}B_k\alpha^{m+1-k} + (m+1)\overline{B}_1\alpha^m - (m+1)\alpha^m \right\}$

$$= -\tfrac{1}{m+1}\left\{\sum_{\substack{k=0\\k\neq 1}}^{m+1}\binom{m+1}{k}B_k\alpha^{m+1-k} + (m+1)(1+B_1)\alpha^m - (m+1)\alpha^m\right\} = -\tfrac{1}{m+1}\left\{\sum_{k=0}^{m+1}\binom{m+1}{k}B_k\alpha^{m+1-k}\right\}$$

Thus $-\frac{B_{m+1}(\alpha)}{m+1} = \zeta(-m,\alpha) = -\tfrac{1}{m+1}\left\{\sum_{k=0}^{m+1}\binom{m+1}{k}B_k\alpha^{m+1-k}\right\}$ for $m \geq 0$

so that $B_{m'}(\alpha) = \sum_{k=0}^{m'}\binom{m'}{k}B_k\,\alpha^{m'-k}$ ,where $m' = m+1$ .

**Proof of Theorem** : We have for complex variables $s_1, s_2, s_3$ with Re $s_1$, Re $s_2 \geq 1$ and Re $s_3 > 1$,

$$\zeta(s_1,s_2,s_3) = \sum_{n_1\geq 1} n_1^{-s} \cdot \sum_{n_2>n_1} n_2^{-s_2} \cdot \sum_{n_3>n_2} n_3^{-s_3} \quad .$$

For Re $s_2$, Re $s_3 > 1$,we have $\sum_{n_2>n_1} n_2^{-s_2} \cdot \sum_{n_3>n_2} n_3^{-s_3} = \sum_{n_2>n_1} n_2^{-s_2}\zeta(s_3, n_2+1)$.

If Re $s_2$ is large , consider $\sum_{n_2>n_1} n_2^{-s_2}\zeta(-m_3, n_2+1) = -\sum_{n_2>n_1} n_2^{-s_2} \cdot \frac{\overline{B}_{m_3+1}(n_2)}{m_3+1}$

$$= -\tfrac{1}{m_3+1}\sum_{n_2>n_1} n_2^{-s_2}\sum_{k_3=0}^{m_3+1}\binom{m_3+1}{k_3}\overline{B}_{k_3} n_2^{m_3+1-k_3} = -\tfrac{1}{m_3+1}\sum_{k_3=0}^{m_3+1}\binom{m_3+1}{k_3}\overline{B}_{k_3} \cdot \sum_{n_2>n_1} n_2^{-(s_2-(m_3+1-k_3))}$$

$$= -\tfrac{1}{m_3+1}\sum_{k_3=0}^{m_3+1}\binom{m_3+1}{k_3}\overline{B}_{k_3} \cdot \zeta(s_2-(m_3+1-k_3), n_1+1)\,.$$

Next for Re $s_1 > 1$ , consider $\sum_{n_1\geq 1} n_1^{-s_1}(-\tfrac{1}{m_3+1})\sum_{k_3=0}^{m_3+1}\binom{m_3+1}{k_3}\overline{B}_{k_3}\zeta(-m_2-m_3-1+k_3, n_1+1)\cdot$

$$= (-\tfrac{1}{m_3+1})\sum_{k_3=0}^{m_3+1}\binom{m_3+1}{k_3}\overline{B}_{k_3} \cdot \sum_{n_1\geq 1} n_1^{-s_1}\zeta(-m_2-m_3-1+k_3, n_1+1)\cdot$$

Next for sufficiently large Re $s_1$, consider

$$\sum_{n_1\geq 1} n_1^{-s_1}\cdot\zeta(-m_2-m_3-1+k_3,n_1+1)=-\sum_{n_1\geq 1} n_1^{-s_1}\cdot\frac{\overline{B}_{m_2+m_3+2-k_3}(n_1)}{m_2+m_3+2-k_3}$$

$$=-\sum_{n_1\geq 1} n_1^{-s_1}\cdot\frac{\sum_{k_2=0}^{m_2+m_3+2-k_3}\binom{m_2+m_3+2-k_3}{k_2}\overline{B}_{k_2}\, n_1^{m_2+m_3+2-k_3-k_2}}{(m_2+m_3+2-k_3)}$$

$$=-\sum_{k_2=0}^{m_2+m_3+2-k_3}\binom{m_2+m_3+2-k_3}{k_2}\frac{\overline{B}_{k_2}}{(m_2+m_3+2-k_3)}\cdot\sum_{n_1\geq 1} n_1^{-(s_1-(m_2+m_3+2-k_3-k_2))}$$

$$=-\sum_{k_2}\binom{m_2+m_3+2-k_3}{k_2}\frac{\overline{B}_{k_2}\,\zeta(s_1-(m_2+m_3+2-k_2-k_3)}{m_2+m_3+2-k_3}\ .$$

However , this is valid for any $s_1$ with $s_1-(m_2+m_3+2-k_2-k_3)\neq 1$ .

Thus , if $s_1$ is any complex number with Re $s_1\leq 0$ ; or if Re $s_1>0$ with $s_1$ not an integer , then

$$\zeta(s_1,-m_2,-m_3)=\frac{1}{m_3+1}\sum_{k_3=0}^{m_3+1}\binom{m_3+1}{k_3}\overline{B}_{k_3}\cdot\sum_{k_2=0}^{m_2+m_3+2-k_3}\binom{m_2+m_3+2-k_3}{k_2}\frac{\overline{B}_{k_2}}{m_2+m_3+2-k_3}\cdot$$

$$\cdot\,\zeta(s_1-(m_2+m_3+2-k_2-k_3))\ .$$

This completes the proof of Theorem .